\documentclass[a4paper,12pt,reqno]{amsart}
\usepackage{amsmath}
\usepackage{amssymb}
\usepackage{amsthm}

\usepackage{amscd}
\usepackage{amsfonts}
\usepackage{setspace}
\usepackage{version}

\title[Limit distributions of a random multiplicative function]{On the limit distributions of some sums of a random multiplicative function}
\author{Adam J Harper}
\address{Department of Pure Mathematics and Mathematical Statistics, Wilberforce Road, Cambridge CB\textup{3} \textup{0}WA, England}
\email{A.J.Harper@dpmms.cam.ac.uk}
\date{30th November 2010}
\thanks{The author is supported by a studentship from the Engineering and Physical Sciences Research Council of the United Kingdom.}

\numberwithin{equation}{section}
\theoremstyle{plain}

\onehalfspace
\newcommand{\N}{\mathbb{N}}
\newcommand{\R}{\mathbb{R}}
\newcommand{\E}{\mathbb{E}}
\newcommand{\p}{\mathbb{P}}
\newcommand{\Z}{\mathbb{Z}}

\newtheorem{result}{Theorem}
\newtheorem{resultprop}{Proposition}
\newtheorem{result2}[result]{Theorem}
\newtheorem{resultcor}{Corollary}
\newtheorem{result3}[result]{Theorem}
\newtheorem{mclt}{Central Limit Theorem}
\newtheorem{wbound}{Number Theory Result}
\newtheorem{wbound2}[wbound]{Number Theory Result}
\newtheorem{wbound3}[wbound]{Number Theory Result}
\newtheorem{wbound4}[wbound]{Number Theory Result}
\newtheorem{wbound5}[wbound]{Number Theory Result}
\newtheorem{techlemma}{Technical Lemma}
\newtheorem{lemma}{Lemma}

\begin{document}

\maketitle

\begin{abstract}
We study sums of a random multiplicative function; this is an example, of number-theoretic interest, of sums of products of independent random variables (chaoses). Using martingale methods, we establish a normal approximation for the sum over those $n \leq x$ with $k$ distinct prime factors, provided that $k=o(\log\log x)$ as $x \rightarrow \infty$. We estimate the fourth moments of these sums, and use a conditioning argument to show that if $k$ is of the order of magnitude of $\log\log x$ then the analogous normal limit theorem does not hold. The methods extend to treat the sum over those $n \leq x$ with at most $k$ distinct prime factors, and in particular the sum over all $n \leq x$. We also treat a substantially generalised notion of random multiplicative function.
\end{abstract}

\section{Introduction}
Let $\epsilon_{p}$ be a sequence of independent Rademacher random variables, indexed by primes $p$; that is, let
$$\p(\epsilon_{p}=1)=\p(\epsilon_{p}=-1)=1/2,$$
independently for each $p$. We construct a random multiplicative function $f$, in the sense of Hal\'{a}sz~\cite{hal} and Wintner~\cite{wint}, by defining
$$f(n) := \left\{ \begin{array}{ll}
     \prod_{p|n} \epsilon_{p} & \textrm{if } n \textrm{ is squarefree}   \\
     0 & \textrm{otherwise}
\end{array} \right.$$

In 1944, Wintner~\cite{wint} studied the behaviour of the summatory function $M(x):=\sum_{n \leq x} f(n)$, as a heuristic for the behaviour of Mertens' function\footnote{The Mertens function is the summatory function of the M\"{o}bius function $\mu(n)$, the multiplicative function taking value $-1$ on all primes, and supported on squarefree numbers only. As discussed, for example, in Chapter 14 of Titchmarsh~\cite{tm}, showing that $\sum_{n \leq x} \mu(n) = O(x^{1/2+\epsilon})$ for each $\epsilon > 0$ is equivalent to proving the Riemann Hypothesis.}. He showed, amongst other things, that for each $\epsilon > 0$ one almost surely has
$$M(x)=O(x^{1/2+\epsilon}) \; \textrm{ as } x \rightarrow \infty.$$
This bound was improved by later authors, (see the discussion of Erd\H{o}s in his unsolved problem papers~\cite{erd1,erd2}), and the best currently known appears to be the 1982 result of Hal\'{a}sz~\cite{hal}: there is an absolute and effective constant $A > 0$ such that, almost surely,
$$M(x)=O(\sqrt{x} e^{A\sqrt{\log\log x \log\log\log x}}) \; \textrm{ as } x \rightarrow \infty.$$
However, by this point the motivation for the problem had shifted somewhat, with Hal\'{a}sz~\cite{hal} writing that ``A deeper aspect... is to find out how the number-theoretic dependence among $f(n)$ effects the magnitude, compared especially with the case of $f(n)=\pm 1$ being independent for all $n$''. This is also the attitude that we adopt.

\vspace{12pt}
Whereas Hal\'{a}sz~\cite{hal} was making a comparison with the law of the iterated logarithm, we will be interested in distributional properties of sums of $f(n)$, comparing with the central limit theorem. Let $\omega(n)$ denote the number of distinct prime factors of $n$, (so e.g. $\omega(12)=2$), and for each $k \in \N$ define
$$M^{(k)}(x):= \sum_{n \leq x, \omega(n)=k} f(n), \;\; \textrm{ and } \; \; \widetilde{M}^{(k)}(x) := M^{(k)}(x)/\sqrt{\E M^{(k)}(x)^{2}}.$$
For these sums we have the following positive result. 
\begin{result}
Suppose that $k \geq 1$ is $o(\log\log x)$ as $x \rightarrow \infty$. Then for each $z \in \R$,
$$\p \left( \widetilde{M}^{(k)}(x) \leq z \right) \rightarrow \Phi(z)$$
as $x \rightarrow \infty$, where $\Phi(z) = \frac{1}{\sqrt{2\pi}} \int_{-\infty}^{z} e^{-t^{2}/2} dt$ is the standard normal distribution function.
\end{result}
Theorem 1 is a refinement of a result of Hough~\cite{ho,ho2}, who showed by the method of moments that taking $k \geq 1$, $k=o(\log\log\log x)$ is permissible for a normal approximation. The theorem will be proved in $\S 4$, and we postpone until then any motivation for studying subsums of this type (but see below).

Readers may identify $M^{(k)}(x)$ as an example of so-called {\em Rademacher chaos of order $k$}, or as a generalised type of {\em U-statistic}\footnote{A $U$-statistic of order $k$ has the form $\frac{(n-k)!}{n!} \sum h(X_{i(1)},...,X_{i(k)})$, where $(X_{i})$ is an i.i.d. sequence of random variables, $h$ is a real-valued function, and the sum is over all tuples $(i(1),...,i(k))$ of distinct numbers smaller than $n$. See the book of de la Pe\~{n}a and Gin\'{e}~\cite{dlpg} for much more discussion of $U$-statistics and chaoses.}. Our proof of Theorem 1 uses a martingale version of the central limit theorem, due to McLeish~\cite{mcl}, which is an idea obtained by the author after reading the paper of Blei and Janson~\cite{bj}. They apply martingale methods to study a general Rademacher chaos, but, in common with other articles on Rademacher chaos, do this when the order $k$ is fixed. In our special case, we apply information about numbers with a given quantity of prime factors, and find we can let $k(x)$ tend to infinity along with $x$.

In the framework of martingale theory, the computations that allow us to deduce Theorem 1 imply other results about $M^{(k)}(x)$. For example, combining the work of $\S 4$ with a central limit theorem of Haeusler and Joos~\cite{hj}, we obtain (roughly) a rate of convergence $O\left((1+z^{4})^{-1} (k/\log\log x)^{2/5} \right)$ in Theorem 1. We refer the reader to Hall and Heyde's book~\cite{hh} for further discussion of the sorts of result that are possible.

\vspace{12pt}
It is a classical result of Hardy and Ramanujan~\cite{hr} that ``the normal order of the number of different prime factors of a number is $\log\log n$'', and in particular
$$\#\{n \leq x : |\omega(n)-\log\log x| > (\log\log x)^{3/4}\} = o(x) \; \textrm{ as } x \rightarrow \infty.$$
This means that the range of $k$ allowed in Theorem 1 stops just short of telling us about sums of $f(n)$ over `typical' numbers.

As we will discuss in $\S 5$, there is a clear change of behaviour of $M^{(k)}(x)$ when $k$ is comparable in size with $\log\log x$, compared to when $k=o(\log\log x)$. One aspect is that, if $p$ is some `not very large' prime, one has
$$\#\{n \leq x: \omega(n)=k(x), p|n\} = o\left( \frac{1}{p} \#\{n \leq x: \omega(n)=k(x)\} \right)$$
if $k(x)=o(\log\log x)$, whereas numbers with about $\log\log x$ prime factors are divisible by such $p$ in roughly the usual proportions $1/p$. Thus, in the latter case, there are many $\epsilon_{p}$ with a large influence on the behaviour of $M^{(k)}(x)$. This is evidenced by the following moment estimate: the reader should observe that the quantity $L$ appearing is of the order of magnitude of $\log\log x$, so the quantity $\overline{k}$ is of the order of magnitude of $k/\log\log x$, whenever $k \leq \log^{0.9}x$, say.
\begin{resultprop}
Let $k \geq 2$, $x \geq 3$, and write
$$m_{4}^{(k)}(x) := \E\widetilde{M}^{(k)}(x)^{4} - 3.$$
Also write $L=L(k,x):=\log\log x - \log k - \log\log (k+1)$, and $\overline{k}:=k/L$. There are constants $C, \delta > 0$ such that, provided $x \geq C$:
\begin{enumerate}
\item if $k \leq (\log\log x)^{1.9}$, then $m_{4}^{(k)}(x) \gg (\overline{k}/\log (\overline{k}+2))^{2}$;
\item if $k \leq \delta \log x/(\log\log x)^{2}$, then $m_{4}^{(k)}(x) \gg (\overline{k}/\log k)^{2}$.
\end{enumerate}
\end{resultprop}
The first bound is stronger than the second roughly when $k \leq (\log\log x)^{1+o(1)}$, and otherwise the second bound extends the first. If $\epsilon > 0$, and $k \leq (1/2-\epsilon)\log\log x$, one can use the methods of $\S 4$ to obtain upper bounds for $m_{4}^{(k)}(x)$ that agree, up to an $\epsilon$-dependent constant, with these lower bounds. For larger $k$ it seems unlikely that Proposition 1 is sharp.

\vspace{12pt}
The inflation of the fourth moment of $\widetilde{M}^{(k)}(x)$, compared with $\E N(0,1)^{4}=3$, suggests that the range of $k$ in Theorem 1 may be best possible. However, it is easy to construct examples of random variables for which this intuition fails. In $\S 6$, we will prove that this is not the case for us.
\begin{result2}
Let $\epsilon, R$ be fixed such that $0 < \epsilon < R$, and suppose that for all large $x$,
$$\epsilon \log\log x \leq k(x) \leq R \log\log x.$$
Then Theorem 1 does {\em not} hold for $\widetilde{M}^{(k)}(x)$ as $x \rightarrow \infty$.
\end{result2}
The fact that the fourth moment of $\widetilde{M}^{(k)}(x)$ does not converge to 3 would establish Theorem 2, if some higher absolute moment was bounded as $x \rightarrow \infty$. See Chapter 8.1 of Feller's book~\cite{fel}. However, the author did not explore this approach, for two reasons: firstly, because such boundedness would very likely not hold on the whole range of $k$ in Theorem 2 (by analogy with the situation for $M(x)$), and would certainly be difficult to establish. In general, it is interesting to consider how to disprove convergence in distribution when computing even moments does not help. Our proof of Theorem 2 uses a conditioning argument, which supplies upper bounds for some truncated moments of $\widetilde{M}^{(k)}(x)$. The inflation of $\E \widetilde{M}^{(k)}(x)^{4}$ means that the distribution of $\widetilde{M}^{(k)}(x)$ has greater {\em kurtosis} than the standard normal distribution--- that is, sharper concentration of most of the density around the mean, and heavier tails. The fact underlying Theorem 2 is that, in this case, the tails are heavy enough to noticeably depress some truncated averages.

\vspace{12pt}
The preceding arguments extend fairly easily to deal with some other types of sum. Introduce the notation
$$M^{(\leq k)}(x):= \sum_{n \leq x, \omega(n) \leq k} f(n), \;\; \textrm{ and } \;\; \widetilde{M}^{(\leq k)}(x):= M^{(\leq k)}(x)/\sqrt{\E M^{(\leq k)}(x)^{2}}.$$
Ben Green suggested these sums to the author as an object of study, and for them we have the following result.
\begin{resultcor}
If $\widetilde{M}^{(k)}(x)$ is replaced by $\widetilde{M}^{(\leq k)}(x)$, then:
\begin{enumerate}
\item Theorem 1 holds without change;

\item Theorem 2 holds without the need for an upper bound on $k$, i.e. it is enough if
$$k(x) \geq \epsilon \log\log x$$
for fixed $\epsilon > 0$. 

\end{enumerate}
\end{resultcor}
Notice that the complete sum $M(x)$ is the same thing as $M^{(\leq x)}(x)$, so in particular $M(x)/\sqrt{\E M(x)^{2}}$ does not converge in distribution to $N(0,1)$. This confirms a heuristic of Chatterjee, expressed in Hough's paper~\cite{ho} (and also see the very recent preprint~\cite{chatsound} of Chatterjee and Soundararajan). See $\S 6$ for some more discussion of this.

A few remarks should be enough to convince the reader that Corollary 1 is true. If $k=o(\log\log x)$, then the number of integers with at most $k-1$ prime factors is of smaller order than the number with exactly $k$ prime factors. This is basically enough for the extension of Theorem 1, whilst extending Theorem 2 is accomplished by summing bounds obtained in the original proof. We provide a few more details of these arguments in $\S 7$. One could extend Proposition 1 in a similar way, but there seems to be little interest in doing this, as the lower bounds obtained would be rather weak.

\vspace{12pt}
The reader may notice, when reading $\S\S \textrm{4--7}$, that the fact that $\epsilon_{p}$ are Rademacher random variables is not used in a very essential way. Our final result establishes a precise version of this principle.
\begin{result3}
Let $\epsilon_{p}$ be any sequence of {\em independent} random variables, also satisfying:
\begin{enumerate}
\item {\em (symmetry)} $\epsilon_{p}$ has the same distribution as $-\epsilon_{p}$, for each $p$;
\item {\em (normalisation)} $\E(\epsilon_{p}^{2})=1$ for each $p$;
\item {\em (bounded fourth moments)} there is a constant $C>0$ such that $\E(\epsilon_{p}^{4}) \leq C$ for all $p$.
\end{enumerate}
If a random multiplicative function is defined using the $\epsilon_{p}$, then the results of Theorem 1, Proposition 1, Theorem 2, and Corollary 1 still hold.
\end{result3}
The fact that Proposition 1 continues to hold is immediate, because $\E(\epsilon_{p}^{4}) \geq (\E(\epsilon_{p}^{2}))^{2}=1$, so the case of Rademacher random variables is worst possible for obtaining lower bounds on fourth moments. Obtaining Theorems 1 and 2 in general requires more work, and we sketch the necessary modifications in $\S 8$. We leave it to the reader to verify that Corollary 1 holds in the more general setting.

One important case included in Theorem 3 is that of standard normal random variables $\epsilon_{p}$, leading to a random multiplicative function that is an example of Gaussian chaos, but we emphasise that there is no requirement even that the $\epsilon_{p}$ be identically distributed.

\vspace{12pt}
We conclude this introduction with a brief list of further work suggested by our results. If we combine Theorems 1 and 2, we still have no information about the behaviour of $M^{(k)}(x)$ if $k(x)$ is not bounded by a constant multiple of $\log\log x$. Although very few numbers $n \leq x$ have so many distinct prime factors, it would be nice to have a complete result. It seems likely, in view of Proposition 1, that the negative result in Theorem 2 could be extended onto the larger range.

Another problem is to give a positive description of the limiting behaviour of $M^{(k)}(x)$ when $\epsilon \log\log x \leq k(x) \leq R \log\log x$. It may be unreasonable to expect a simple limit distribution, and the type of theorem sought might say that one obtains the same limiting distribution for any $\epsilon_{p}$ as in Theorem 3. A reader interested in pursuing this might consult, for example, the book of de la Pe\~{n}a and Gin\'{e}~\cite{dlpg}, who present similar theorems for $U$-statistics.

A third line of work would be to return to the problem treated by Hal\'{a}sz~\cite{hal} and others, of obtaining almost sure bounds for $M(x)$. As remarked previously, (and also see $\S\S 2,6$), there are many random variables $\epsilon_{p}$ having a large influence on this sum, so it does not have the `concentration' properties under which martingale methods perform nicely. Because of this, obtaining sharp bounds for $M(x)$ seems to remain difficult. 

\section{The martingale central limit theorem}
We begin by applying some observations of Blei and Janson~\cite{bj}, made by them in the context of a general Rademacher chaos of fixed order. Let $P(n)$ denote the largest prime factor of $n \geq 2$, and for prime $p$ write
$$M_{p}^{(k)}(x) := \sum_{n \leq x, \omega(n)=k, \atop P(n)=p} f(n).$$
Also let $(\mathcal{F}_{p})$ denote the natural filtration of the random variables $\epsilon_{p}$, so that $\mathcal{F}_{p}$ is the sigma algebra generated by $\{ \epsilon_{q} : q \leq p\}$.

We have the decomposition
$$M^{(k)}(x) = \sum_{p \leq x} M_{p}^{(k)}(x).$$
Since $\epsilon_{p}$ occurs in each summand of $M_{p}^{(k)}(x)$, and not in $M_{q}^{(k)}(x)$ for any $q < p$, the conditional expectation
$$\E\left(M_{p}^{(k)}(x) \mid \{\epsilon_{q} : q < p\}\right)$$
vanishes, so $\left(M_{p}^{(k)}(x) \right)$ is a {\em martingale difference sequence} relative to $((\mathcal{F}_{p}), \p)$.

Like Blei and Janson~\cite{bj}, we will use the following version, essentially Corollary 2.13 of McLeish~\cite{mcl}, of the martingale central limit theorem.
\begin{mclt}[McLeish, 1974]
For $n \in \N$, suppose that $k_{n} \in \N$, and that $X_{i,n}$, $1 \leq i \leq k_{n}$ is a martingale difference sequence on $(\Omega,\mathcal{F},(\mathcal{F}_{i,n})_{i},\p)$. Write $S_{n}:=\sum_{i \leq k_{n}} X_{i,n}$, and suppose that the following conditions hold:
\begin{enumerate}
\item {\em (normalised variances)} $\sum_{i \leq k_{n}} \E X_{i,n}^{2} \rightarrow 1$ as $n \rightarrow \infty$;
\item {\em (Lindeberg condition)} for each $\epsilon > 0$, we have $\sum_{i \leq k_{n}} \E(X_{i,n}^{2} \textbf{1}_{|X_{i,n}|>\epsilon}) \rightarrow 0$ as $n \rightarrow \infty$;
\item {\em (cross terms condition)} $\limsup_{n \rightarrow \infty} \sum_{i \leq k_{n}} \sum_{j \leq k_{n}, j \neq i} \E X_{i,n}^{2} X_{j,n}^{2} \leq 1$.
\end{enumerate}
Then $S_{n}$ converges in distribution to $N(0,1)$ as $n \rightarrow \infty$.
\end{mclt}
We will apply this to the normalised random variables
$$M_{p}^{(k)}(x) / \sqrt{\E M^{(k)}(x)^{2}},$$
where $x$ (which we can obviously allow to tend to infinity through integer values) takes the place of $n$, and primes $p$ play the role of the indices $i$. In fact, since
\begin{eqnarray}
\E M^{(k)}(x)^{2} & = & \sum_{m \leq x, \omega(m) = k} \sum_{n \leq x, \omega(n) = k} \E f(m)f(n) \nonumber \\
& = & \#\{n \leq x : n \textrm{ is squarefree}, \omega(n)=k\}, \nonumber
\end{eqnarray}
and similarly for $\E M_{p}^{(k)}(x)^{2}$, the normalised variances condition holds trivially (i.e. not just in the limit, but for each $x$ individually\footnote{Here we have only used that the random variables $f(n)$ are {\em orthogonal}, i.e. $\E f(m)f(n) = 0$ unless $m=n$.}).

We comment briefly on the conditions of Central Limit Theorem 1, although the reader is referred again to McLeish's excellent paper~\cite{mcl} for further discussion. The Lindeberg condition forces that individual summands $X_{i,n}$ are ``asymptotically negligible''; the collective force of the conditions is, roughly, that
$$\E \left(\sum_{i \leq k_{n}} X_{i,n}^{2} - 1 \right)^{2} \rightarrow 0 \textrm{ as } n \rightarrow \infty,$$
so that the sum of squares concentrates around the desired variance in the limit. The failure of this behaviour will feature prominently in $\S 6$.

\section{Some number-theoretic estimates}
Next, we record some estimates for the quantity of (squarefree) numbers with a fixed number of distinct prime factors, which will be needed repeatedly. We write $\Omega(n)$ for the number of prime factors of $n$ counted with multiplicity, so e.g. $\Omega(12)=3$. Thus $\omega(n)=\Omega(n)$ if and only if $n$ is squarefree, and this is how the notation $\Omega(\cdot)$ will be used.

The following is a standard upper bound:
\begin{wbound}[Hardy and Ramanujan, 1917]
There are absolute constants $A, B$ such that, for all $k \geq 1$ and $x \geq 2$,
$$\#\{n \leq x : \omega(n)=k\} \leq \frac{Ax(\log\log x + B)^{k-1}}{(k-1)!(\log x)}.$$
\end{wbound}
This is due to Hardy and Ramanujan~\cite{hr}, and the reader may wish to note (although we do not need to do so here) that analogous uniform bounds do not hold when counting prime factors with multiplicity, whilst in many other situations there is essentially no difference between these cases.

We cite the following lower bound from the text of Montgomery and Vaughan~\cite{mv}; it is easily implied, for example, by Exercise 4 of their $\S 7.4$.
\begin{wbound2}[Sathe, Selberg, 1954]
There is a constant $\delta > 0$ such that, for $x \geq 3$ and $1 \leq k \leq \log\log x$,
$$\#\{n \leq x : \omega(n)=\Omega(n)=k\} \geq \delta \frac{x (\log\log x)^{k-1}}{(k-1)! (\log x)}.$$
\end{wbound2}
Montgomery and Vaughan~\cite{mv} present Selberg's proof, refining one of Sathe, of asymptotic formulae for the quantities in Number Theory Results 1 and 2. This is more than we need, but the author knows of no earlier proof obtaining lower bounds on a comparable range of $k$.

\vspace{12pt}
In $\S\S \textrm{5--8}$, we will need some variants of these two results. We begin by presenting an `elementary' lower bound of Pomerance~\cite{pom}:
\begin{wbound3}[Pomerance, 1984]
There is an absolute constant $C$ such that, if $x \geq C$ and $\log\log x (\log\log\log x)^{2} \leq k \leq \log x /(3 \log\log x)$, we have
$$\#\{n \leq x: \omega(n)=\Omega(n)=k\} \geq \frac{x}{k! \log x} e^{k(\log L + \log L /L +O(1/L))},$$
where $L=L(k,x):=\log\log x - \log k - \log\log (k+1)$.
\end{wbound3}
We will use this to verify that $\#\{n \leq x: \omega(n)=\Omega(n)=k\}$ is much larger than some error terms that will be subtracted from it. Observe that the upper bound restriction on $k$ is impressively large, differing by a bounded factor from the maximum for which the left hand side is non-zero.

In their 1988 paper, Hildebrand and Tenenbaum~\cite{ht} give various estimates for $\#\{n \leq x : \omega(n)=k\}$ on the range $1 \leq k \leq \log x /(\log\log x)^{2}$, which appear to include all previous results on that range. We will make substantial use of slight modifications of these; but it is fairly easy to see that the methods of Hildebrand and Tenenbaum~\cite{ht} imply the results. (Most of their arguments carry over without change, once one checks that suitable parts of their Lemmas 1 and 2 hold in the modified situations).

Hildebrand and Tenenbaum~\cite{ht} use their asymptotic for $\#\{n \leq x : \omega(n)=k\}$, which involves two `implicit' parameters $\rho=\rho(x,k)$ and $\alpha=\alpha(x,k)$, to ``...get rather precise information on the local behaviour of the function''. In an exactly similar way, we have:
\begin{wbound4}[Hildebrand and Tenenbaum, 1988]
There are absolute constants $C, \delta > 0$ such that, if $x \geq C$ and $1 \leq k \leq \delta \log x /(\log\log x)^{2}$,
$$\frac{\#\{n \leq x: \omega(n)=\Omega(n)=k+1\}}{\#\{n \leq x: \omega(n)=\Omega(n)=k\}} = \frac{L}{k}\left(1+O\left(\frac{\log L}{L} \right)\right),$$
where $L$ is as in Number Theory Result 3. If $1 \leq \lambda \leq x$, also
$$\frac{\#\{n \leq \lambda x: \omega(n)=\Omega(n)=k\}}{\#\{n \leq x: \omega(n)=\Omega(n)=k\}} = \lambda \left(1+\frac{\log \lambda}{\log x}\right)^{(k/L)-1} e^{O(1/L + k \log L \log \lambda/(L^{2}\log x))}.$$
\end{wbound4}

\vspace{12pt}
Finally, we state a rather specialised variant of Corollary 1 of Hildebrand and Tenenbaum~\cite{ht}; it is possible to give an explicit expression for the function $h=h(k,x)$ that appears, but we simply note that $h=O(\log^{2}(2+\overline{k})/(\log\log x)^{2})$.
\begin{wbound5}
Let $R \geq 0$, and let $\mathcal{P}$ be a set of prime numbers. There are absolute constants $C, \delta > 0$ such that, if $x \geq C$ and $1 \leq k \leq \delta (\log\log x)^{2}$,
$$ \#\{n \leq x: \omega(n)=\Omega(n)=k, \;\; p \nmid n \; \forall p \in \mathcal{P} \}$$
$$ = G(\overline{k}) \prod_{p \in \mathcal{P}} \left(1+\frac{\overline{k}}{p} \right)^{-1} \frac{x (\log\log x)^{k-1}}{(k-1)! \log x} e^{-kh/2} \left(1+O\left(\frac{1}{\log\log x} \right) + O\left(\frac{k}{(\log\log x)^{2}} \right) \right),$$
uniformly for $\#\mathcal{P} \leq R$. Here $\overline{k}$ is as in the Introduction, and
$$G(z):=\frac{1}{\Gamma(z+1)} \prod_{p \textrm{ prime}} \left(1+\frac{z}{p} \right) \left(1-\frac{1}{p} \right)^{z}.$$
\end{wbound5}
When thinking about the function $G(z)$, the following estimate will occasionally be useful\footnote{This is modified, and slightly corrected, from p 498 of Hildebrand and Tenenbaum's paper~\cite{ht}.}: uniformly for $z \geq 2$, we have
\begin{eqnarray}
\frac{d}{dz} \log(zG(z)) & = & -\frac{\Gamma'(z)}{\Gamma(z)} + \sum_{p \textrm{ prime}} \left(\frac{1}{p+z} + \log \left(1-\frac{1}{p} \right) \right) \nonumber \\
& = & -\log(z\log z) - \gamma + O\left( \frac{1}{\log z} \right), \nonumber
\end{eqnarray}
where $\gamma \approx 0.577$ is Euler's constant.

\section{Proof of Theorem 1}

\subsection{Some opening remarks}
Before embarking on the proof of Theorem 1, we give the promised motivation for studying the sums $M^{(k)}(x)$. One reason is the connection with the large body of probabilistic literature on Rademacher chaos: it is usual to study chaoses of fixed order $k$, cf. the papers of Blei and Janson~\cite{bj} or of Nourdin, Peccati and Reinert~\cite{npr}, although one does not typically let $k=k(x)$ tend to infinity. When Hough~\cite{ho} studies $M^{(k)}(x)$, he wishes to exploit that ``...for $k$ slowly growing compared to $x$, numbers having only $k$ prime factors almost always have a [very large] prime factor...''. However, there is a third justification for studying $M^{(k)}(x)$, which appeals to the present author rather more.

Fix $x$ and $k$, and let $I$ be a random variable independent of all the $\epsilon_{p}$, having the discrete uniform distribution on
$$\{p \leq x : p \textrm{ is prime}\}.$$
Define a new random multiplicative function $f'(n)$, by replacing $\epsilon_{p}$ with an independent copy $\epsilon_{p}'$ when $I=p$, and leaving the other $\epsilon_{q}$ unchanged. If we write $N^{(k)}(x) := \sum_{n \leq x, \omega(n)=k} f'(n)$, then the pair $(M^{(k)}(x),N^{(k)}(x))$ has two important properties:
\begin{enumerate}
\item (exchangeability) for all sets $C,D \subseteq \Z$,
$$\p(M^{(k)}(x) \in C \textrm{ and } N^{(k)}(x) \in D) = \p(M^{(k)}(x) \in D \textrm{ and } N^{(k)}(x) \in C);$$

\item (regression) writing (as usual) $\pi(x) := \#\{p \leq x : p \textrm{ is prime}\}$,
$$\E(N^{(k)}(x) | M^{(k)}(x)) = \left(1-\frac{k}{\pi(x)} \right) M^{(k)}(x).$$

\end{enumerate}
We leave the reader to verify these, and that an analogue of the regression property would {\em not} hold if one attempted a similar construction for e.g. the complete sum $M(x)$. We will not use the exchangeable pair construction here, but it is often useful to have it available: see, for example, the application to normal approximation in $\S 3.3$ of Nourdin, Peccati and Reinert's paper~\cite{npr}.

\vspace{12pt}
The proof of Theorem 1 is somewhat technical, so first we highlight what the important steps will be. In $\S 4.2$, we will reduce the proof to showing that a certain sum, counting quadruples of integers according to their `squarefree parts', is of small order when $k(x)=o(\log\log x)$.

We will split this sum, fixing the number of prime factors of the squarefree parts, and in $\S 4.3$ will bound the subsums. The details are rather heavy, but there are many symmetries inherent in our counting procedure, and we try to give a conceptual explanation of this. The results of $\S 4.3$ will be valid for arbitrary $k(x)$.

In $\S 4.4$, we will finish the proof just by summing the upper bounds previously obtained, and using the assumption $k(x)=o(\log\log x)$. At the beginning of $\S 5$, we give a heuristic justification for why (the proof of) Theorem 1 works out, which we hope the reader will also find helpful.

\subsection{Preliminaries to the proof}
Recall that, to deduce Theorem 1 from Central Limit Theorem 1, it remains to verify the Lindeberg and cross terms conditions. As noticed in general by Blei and Janson~\cite{bj}, for any $\epsilon > 0$ we have
$$\sum_{p \leq x} \E\left(\frac{M_{p}^{(k)}(x)^{2}}{\E M^{(k)}(x)^{2}} \textbf{1}_{|M_{p}^{(k)}(x)|/\sqrt{\E M^{(k)}(x)^{2}}> \epsilon} \right) \leq \frac{1}{\epsilon^{2}} \sum_{p \leq x} \E\left(\frac{M_{p}^{(k)}(x)^{4}}{(\E M^{(k)}(x)^{2} )^{2}} \right).$$
Thus the Lindeberg condition certainly holds if, as $x \rightarrow \infty$,
$$\sum_{p \leq x} \E M_{p}^{(k)}(x)^{4} = o(\#\{n \leq x : \omega(n)=\Omega(n)=k\}^{2}).$$

In order to verify the cross terms condition, we must also study
$$\E M_{p}^{(k)}(x)^{2} M_{q}^{(k)}(x)^{2} = \sum_{a \leq x, \omega(a)=k, \atop P(a)=p} \sum_{b \leq x, \omega(b)=k, \atop P(b)=p} \E f(a)f(b) M_{q}^{(k)}(x)^{2} ,$$
where $p$ and $q$ are primes. Introduce the notation
$$S_{p,k,x}:=\{n \leq x : \omega(n)=\Omega(n)=k, P(n)=p\},$$
and write $S_{k,x}$ for the larger set without any restriction on $P(n)$. Notice that we can write $\E M_{p}^{(k)}(x)^{2} M_{q}^{(k)}(x)^{2}$ as
$$\sum_{m \leq x^{2}, \omega(m)=\Omega(m)} \#\{(a,b) \in S_{p,k,x}^{2} : s(ab)=m\} \#\{(c,d) \in S_{q,k,x}^{2} : s(cd)=m\}, $$
where $s(x)$ denotes $x$ divided by its largest square factor, e.g. $s(120)=30$.

The $m=1$ term in the sum is $\#S_{p,k,x} \cdot \#S_{q,k,x}$; summing this over all $p$ and $q$ yields $(\#S_{k,x})^{2}$, whilst for $k=o(\log\log x)$ we have (if $k \geq 2$)
\begin{eqnarray}
\sum_{p \leq x} (\#S_{p,k,x})^{2} \leq \sum_{p \leq x} (\#S_{k-1,x/p})^{2} & = & O\left(\frac{x^{2} (\log\log x + B)^{2k-4}}{(k-2)!^{2}} \sum_{p \leq x} \frac{1}{p^{2} \log^{2}(x/p+1)} \right) \nonumber \\
& = & O\left(\frac{x^{2} (\log\log x + B)^{2k-4}}{(k-2)!^{2} \log^{2}x} \right) \nonumber \\
& = & o((\#S_{k,x})^{2}), \nonumber
\end{eqnarray}
by Number Theory Results 1 and 2. Thus we will have established the Lindeberg and cross terms conditions \emph{if we show that the $m=1$ term gives the main contribution}, i.e. that
$$\sum_{p,q \leq x} \sum_{W=1}^{k-1} \sum_{m \leq x^{2}, \omega(m)=\Omega(m)=2W} \#\{(a,b) \in S_{p,k,x}^{2} : s(ab)=m\} \#\{(c,d) \in S_{q,k,x}^{2} : s(cd)=m\}$$
is $o((\#S_{k,x})^{2})$ as $x \rightarrow \infty$. Values of $m$ with an odd number of prime factors cannot arise as values of $s(ab)$ here: this will be clear when we begin to estimate the sums.

\vspace{12pt}
We must record some technical estimates that will be needed to complete the above programme. These are encapsulated in the following result.
\begin{techlemma}
Let $a \in \N$, $n \in \N \cup \{0\}$, $C \geq 0$, and $2 \leq m \leq M$. Also let $t \in \N \cup \{0\}$, and suppose that $t \leq D \log(N/m)$, where $D \geq 0$ and $N \geq 3m$. Then the following hold:
\begin{enumerate}
\item $\sum_{\sqrt{m} \leq j \leq m, \omega(j)=a} \frac{(\log\log(m/j+2) + C)^{n}}{j \log^{2}(m/j+1)} = O\left(\frac{n!(\log\log m + B)^{a-1}}{(a-1)!\log m} \right)$;
\item $\sum_{j \leq m, \omega(j)=a} \frac{(\log\log(j+1) + C)^{n}}{j \log j} = O\left(\frac{(n+a-1)!}{(a-1)!} \right)$;
\item $\sum_{m \leq j \leq M, \omega(j)=a} \frac{1}{\log^{2}(1+j/m)} = O\left(\frac{M(\log\log M + B)^{a-1}}{(a-1)!\log M \log^{2}(1+M/m)} \right)$;
\item $\sum_{j \leq m, \omega(j)=a} (\log\log(N/j)+C)^{t} = O\left(\frac{m(\log\log m + B)^{a-1}(\log\log(N/m)+C)^{t}}{(a-1)!\log m} \right)$;
\end{enumerate}
where the constants implicit in the ``big Oh'' notation depend at most on $C,D$.
\end{techlemma}

An interested reader will find sketch proofs of these in the appendix.

\subsection{Bounds for the sums with $W$ fixed}
In this section we will prove:
\begin{lemma}
Uniformly for $x \geq 3$, $1 \leq W \leq k-1$; and with the conventions that $(-1)!,(-2)!=1$, and that the third summand is omitted when $j=0$; we have
$$\sum_{p \leq x} \sum_{q \leq p} \sum_{m \leq x^{2}, \omega(m)=\Omega(m)=2W} \#\{(a,b) \in S_{p,k,x}^{2} : s(ab)=m\} \#\{(c,d) \in S_{q,k,x}^{2} : s(cd)=m\}$$
$$ = O\left(\frac{x^{2}}{\log^{2}x (k-W-1)!^{2}} \left( \sum_{j=0}^{W-1} \frac{(2j-2)!(2k-2W-2)!(\log\log x + B)^{2W-2j-2}}{(W-j-1)!^{2} (j-1)!^{2}} + \right. \right.$$
$$ \left. \left. + \frac{(2j-2)!(2W-2j-2)!(\log\log x + B)^{2k-2W-2}}{(W-j-1)!^{2} (j-1)!^{2}} + \right. \right.$$
$$ \left. \left. + \frac{(2W-2j-2)!(2k-2W-2)!(\log\log x + B)^{2j-2}}{(W-j-1)!^{2} (j-1)!^{2}} \right) \right).$$
\end{lemma}

Slightly extending the notation of $\S 4.2$, we define
$$S_{<p,k,x}:=\{n \leq x : \omega(n)=\Omega(n)=k, P(n)<p\},$$
and similarly define $S_{\leq p,k,x}$. Our first important observation is that if squarefree $m$ satisfies $\omega(m)=2W$, then
\begin{eqnarray}
\#\{(a,b) \in S_{p,k,x}^{2} : s(ab)=m\} & = & \#\{(a,b) \in S_{<p,k-1,x/p}^{2} : s(ab)=m\} \nonumber \\
& \leq & \sum_{d|m; d, m/d \leq x/p; \atop \omega(d)=W} \# S_{<p,k-1-W, \min\{x/pd,xd/pm\}}. \nonumber
\end{eqnarray}
The reasoning here justifies our remark that we only need consider $m$ with an even number of prime factors: if $s(ab)=m$ we must have $a=yd$, $b=y(m/d)$ for some $d|m$ and some $y$, (since $a,b$ are squarefree), and since $\omega(a)=\omega(b)$ then also $\omega(d)=\omega(m/d)$.

For the next few lines we will write $S(p,d,m)$ as shorthand for the terms in the sum over $d$; we do not continue to write that we sum over squarefree integers, although we occasionally make use of this fact. We see the sum over $m$ in the statement of Lemma 1 is at most
$$\sum_{d \leq x/p, \atop \omega(d)=W} \sum_{e \leq x/q, \atop \omega(e)=W} \sum_{m \leq xd/p, xe/q; \atop \omega(m)=2W; [d,e]|m} S(p,d,m) \cdot S(q,e,m) $$
$$ \leq \sum_{j=0}^{W} \sum_{h \leq x/p, \atop \omega(h)=j} \sum_{d' \leq x/ph, \atop \omega(d')=W-j} \sum_{e' \leq x/qh, \atop \omega(e')=W-j} \sum_{m' \leq x/pe', x/qd', \atop \omega(m')=j} S(p,hd',hd'e'm') S(q,he',hd'e'm'),$$
where $[d,e]$ denotes the least common multiple of $d,e$, and we reorganise the summations according to the highest common factor $h$ of $d,e$ (putting $d=d'h$, $e=e'h$, and $m=m'de/h$).

The expression above looks difficult to work with, because of the many sums over different ranges. However, by moving the sum over $h$ to the middle we obtain
$$ \sum_{j=0}^{W} \sum_{d' \leq x/p, \atop \omega(d')=W-j} \sum_{e' \leq x/p, \atop \omega(e')=W-j} \sum_{h \leq x/pd', x/qe', \atop \omega(h)=j} \sum_{m' \leq x/pe', x/qd', \atop \omega(m')=j} S(p,hd',hd'e'm') S(q,he',hd'e'm'),$$
where the ranges of summation over $d$ and $e$ are the same. (We can restrict the variable $e'$ to be smaller than $x/p$, rather than $x/q$, because for larger $e'$ the sum over $m'$ is empty). Summing this over $q \leq p$, with the observation that $S(q,he',hd'e'm') = \#S_{q,k-W,\min\{x/he',x/d'm'\}}$, yields
$$ \sum_{j=0}^{W} \sum_{d' \leq x/p, \atop \omega(d')=W-j} \sum_{e' \leq x/p, \atop \omega(e')=W-j} \sum_{h \leq x/pd', x/e', \atop \omega(h)=j} \sum_{m' \leq x/pe', x/d', \atop \omega(m')=j} S(p,hd',hd'e'm') \#S_{\leq p,k-W,\min\{x/he',x/d'm'\}}.$$
(The constraints on $h$ and $m'$ that involve $q$ translate into constraints in the summation over $q$, which are seen to be redundant). Now summing over $p \leq x$, we conclude that the left hand side in Lemma 1 is at most
$$ \sum_{j=0}^{W} \sum_{d' \leq x, \atop \omega(d')=W-j} \sum_{e' \leq x, \atop \omega(e')=W-j} \sum_{h \leq x/d', x/e', \atop \omega(h)=j} \sum_{m' \leq x/e', x/d', \atop \omega(m')=j} \#S_{k-W,\min\{x/hd',x/e'm'\}} \#S_{k-W,\min\{x/he',x/d'm'\}} $$
$$ \leq 4 \sum_{j=0}^{W} \sum_{d' \leq x, \atop \omega(d')=W-j} \sum_{e' \leq d', \atop \omega(e')=W-j} \sum_{h \leq x/d', \atop \omega(h)=j} \#S_{k-W,x/hd'} \left( \sum_{m' \leq e'h/d', \atop \omega(m')=j} \#S_{k-W,x/he'} + \sum_{e'h/d' < m' \leq h, \atop \omega(m')=j} \#S_{k-W,x/d'm'} \right),$$
using the symmetry of $d'$ and $e'$, and of $h$ and $m'$, that emerges.

(To understand the above argument qualitatively, the reader may find it helpful to visualise an annulus, representing a number $m$, split into four sectors: a $j$ prime factor part, representing $h=(d,e)$; two $W-j$ prime factor parts adjoining this, representing $d'=d/h$ and $e'=e/h$; and a $j$ prime factor part representing $m'=m/hd'e'=hm/de$. The symmetry that we have manifested is of the pair $h,d'$ with the pair $m',e'$; for example, although we have $d'h=d \leq x/p$ and $e \leq x/q$, (which may be larger), we also have $e'm'=m/d \leq x/p$.)

\vspace{12pt}
In order to reduce the bound we obtained to the form asserted, we will use Number Theory Result 1 and Technical Lemma 1. It is easiest to treat the terms where $j=0$ and $j=W$ separately; both of these\footnote{The fact that the $j=0$ and $j=W$ contributions are the same, although it can be seen directly, also follows from the ``counting around an annulus'' description. In general, the $j$ and $W-j$ summands in Lemma 1 are the same.} are $\sum_{d' \leq x, \omega(d')=W} (\#S_{k-W,x/d'})^{2} \sum_{e' \leq d', \omega(e')=W} 1$, which is
$$ O \left(\frac{x^{2}}{(k-W-1)!^{2}(W-1)!} \sum_{d' \leq x, \atop \omega(d')=W} \frac{(\log\log(x/d'+1)+B)^{2k-2W-2} (\log\log(d'+1)+B)^{W-1}}{d' \log^{2}(x/d'+1) \log(d'+1)} \right). $$
Splitting the sum at $d'=\sqrt{x}$, and applying the first and second parts of Technical Lemma 1, we obtain the $j=0$ terms in the statement of Lemma 1.

Now fix $1 \leq j \leq W-1$, and consider the second part of the bound, where we sum over large $m'$. This is $x^{2}/(k-W-1)!^{2}$ multiplied by
$$O\left(\sum_{d' \leq x, \atop \omega(d')=W-j} \frac{1}{(d')^{2}} \sum_{h \leq x/d', \atop \omega(h)=j} \textbf{1}_{S_{k-W, x/hd'} \neq \emptyset} \frac{(\log\log(x/hd'+1)+B)^{k-W-1}}{h \log(x/hd'+1)} \right.$$
$$\left. \cdot \sum_{h/d' < m' \leq h, \atop \omega(m')=j} \frac{(\log\log(x/m'd'+1)+B)^{k-W-1}}{m' \log(x/m'd'+1)} \sum_{e' < d'm'/h, \atop \omega(e')=W-j} 1 \right)$$
$$ = O\left(\sum_{d' \leq x, \atop \omega(d')=W-j} \frac{(\log\log(d'+1)+B)^{W-j-1}}{(W-j-1)! d'} \sum_{h \leq x/d', \atop \omega(h)=j} \textbf{1}_{S_{k-W, x/hd'} \neq \emptyset} \frac{(\log\log(x/hd'+1)+B)^{k-W-1}}{h^{2} \log^{2}(x/hd'+1)} \right. $$
$$\left. \cdot \sum_{h/d' < m' \leq h, \atop \omega(m')=j} \frac{(\log\log(x/m'd'+1)+B)^{k-W-1}}{\log(d'm'/h+1)} \right),$$
where in the first bound we move the sum over $e'$ to be performed earlier, because no summand is a function of $e$ except through its range of summation. In the second bound we overestimate $(\log\log(d'm'/h+1)+B)^{W-j-1}$ and $1/\log(x/m'd'+1)$ by their values at the $m'=h$ end of the range of summation. We expect that most of the contribution to the sum over $m'$ will come from this end, since the factors in the summand that decrease with $m'$ do so slowly, and we will see that we do not lose much by doing this.

If $\omega(n)=a$ for some $n \leq y$, then certainly $a \leq \log y/\log 2$, so in estimating the sum over $m'$ we can assume that $k-W \leq \log(x/hd')/\log 2$. Then by the Cauchy--Schwarz inequality, and the third and fourth parts of Technical Lemma 1, the above is
$$ O\left(\frac{1}{(W-j-1)!(j-1)!}\sum_{d' \leq x, \omega(d')=W-j} \frac{(\log\log(d'+1)+B)^{W-j-1}}{d' \log(d'+1)} \right.$$
$$\left. \cdot \sum_{h \leq x/d', \omega(h)=j} \frac{(\log\log(x/hd'+1)+B)^{2k-2W-2}(\log\log(h+1)+B)^{j-1}}{h \log^{2}(x/hd'+1) \log(h+1)} \right), $$
where the value of $B$ is possibly increased by a fixed amount from that in Number Theory Result 1. An exactly similar argument shows that this expression, multiplied by $x^{2}/(k-W-1)!^{2}$, also majorises the first part of our original bound (sum over small $m'$).

Splitting the range of summation over $h$ at $\sqrt{x/d'}$, and applying the first and second parts of Technical Lemma 1 as before, we find this is at most
$$ O\left(\frac{(2j-2)!}{(W-j-1)!(j-1)!^{2}}\sum_{d' \leq x, \omega(d')=W-j} \frac{(\log\log(d'+1)+B)^{W-j-1}}{d' \log(d'+1) \log^{2}(x/d'+1)} (\log\log(x/d'+1)+B)^{2k-2W-2} \right.$$
$$\left. + \frac{(2k-2W-2)!}{(W-j-1)!(j-1)!^{2}}\sum_{d' \leq x, \omega(d')=W-j} \frac{(\log\log(d'+1)+B)^{W-j-1}}{d' \log(d'+1) \log^{2}(x/d'+1)} (\log\log(x/d'+1)+B)^{2j-2} \right).$$
Estimating the sums over $d'$ in the same way proves Lemma 1.
\begin{flushright}
Q.E.D.
\end{flushright}

\subsection{Conclusion of the proof}
In this section we will deduce Theorem 1 from Lemma 1, the preliminary observations in $\S 4.2$, and Number Theory Result 2. In view of these facts, it will suffice to show that the following are $o(1)$ as $x \rightarrow \infty$, when $k(x)=o(\log\log x)$:
$$ \sum_{W=1}^{k-1} \frac{(k-1)!^{2} (2k-2W-2)!}{(k-W-1)!^{2}} \sum_{j=0}^{W-1} \frac{(2j-2)! (\log\log x + B)^{2W-2j-2k}}{(W-j-1)!^{2} (j-1)!^{2}} ;$$
$$ \sum_{W=1}^{k-1} \frac{(k-1)!^{2} (\log\log x + B)^{-2W}}{(k-W-1)!^{2}} \sum_{j=0}^{W-1} \frac{(2j-2)! (2W-2j-2)!}{(W-j-1)!^{2} (j-1)!^{2}} .$$
The third sum in Lemma 1 is bounded above by the first, which can be seen, as remarked earlier, by replacing $j$ by $W-j$ and adding the terms in reverse order.

We illustrate a suitable argument, which is just straightforward analysis, for the first expression. In the inner sum, the ratio of the $j+1$ and $j$ summands is
$$\frac{(2j)! (\log\log x + B)^{2W-2j-2k-2} (W-j-1)!^{2} (j-1)!^{2}}{(2j-2)! (\log\log x + B)^{2W-2j-2k} (W-j-2)!^{2} j!^{2}} \leq \frac{4 (W-j-1)^{2}}{(\log\log x + B)^{2}}.$$
Since $W-j-1 \leq W \leq k$, this is smaller than $1/2$ for $x$ large enough (depending on how quickly $k(x)/\log\log x$ tends to 0 with $x$), and the inner sum is then at most
$$ \frac{2 (\log\log x + B)^{2W-2k}}{(W-1)!^{2}} .$$
Treating the outer sum in the same way, we find it is dominated by the $W=k-1$ term when $x$ is large; and that term is
$$ 2 (k-1)^{2} (\log\log x + B)^{-2} = o(1) \textrm{ as } x \rightarrow \infty.$$

To deal with the second expression, we just note that
$$ \sum_{j=0}^{W-1} \frac{(2j-2)! (2W-2j-2)!}{(W-j-1)!^{2} (j-1)!^{2}} \leq \sum_{j=0}^{W-1} 2^{2W} = W 2^{2W},$$
and then bound the outer sum as before (the $W=1$ term proving to be the dominant one).
\begin{flushright}
Q.E.D.
\end{flushright}

\section{Proof of Proposition 1}
We notice that
\begin{eqnarray}
\E M^{(k)}(x)^{4} = \E(\sum_{p \leq x} M_{p}^{(k)}(x))^{4} & \geq & \sum_{p \leq x} \E M_{p}^{(k)}(x)^{4} + 6 \sum_{p \leq x} \sum_{q < p} \E M_{p}^{(k)}(x)^{2} M_{q}^{(k)}(x)^{2} \nonumber \\
& = & 3\sum_{p \leq x} \sum_{q \leq x} \E M_{p}^{(k)}(x)^{2} M_{q}^{(k)}(x)^{2} - 2\sum_{p \leq x} \E M_{p}^{(k)}(x)^{4}, \nonumber
\end{eqnarray}
since each term in the expansion of $\E(\sum_{p \leq x} M_{p}^{(k)}(x))^{4}$ is non-negative. Referring to the calculations in $\S 4.2$, we see\footnote{One also needs to check that $\sum_{p \leq x} (\#S_{p,k,x})^{2} = o((\#S_{k,x})^{2})$ on the whole range of $k$ in Proposition 1, and not just when $k=o(\log\log x)$. This follows because, for example,
$$ \sum_{p \leq x} (\#S_{p,k,x})^{2} = \sum_{k \leq p \leq x} (\#S_{p,k,x})^{2} \leq \#S_{k-1,x/k} \sum_{k \leq p \leq x} \#S_{p,k,x} = \#S_{k-1,x/k} \cdot \#S_{k,x}, $$
where $\# S_{k-1,x/k} = o(\#S_{k,x})$ by Number Theory Result 4.} that the term $m=1$ there produces the ``main term'' $3$ in Proposition 1, and so to prove Proposition 1 it will suffice to obtain analogous lower bounds for any terms from
$$\sum_{p,q \leq x} \sum_{W=1}^{k-1} \sum_{m \leq x^{2}, \omega(m)=\Omega(m)=2W} \#\{(a,b) \in S_{p,k,x}^{2} : s(ab)=m\} \#\{(c,d) \in S_{q,k,x}^{2} : s(cd)=m\}.$$
Because we are trying to establish lower bounds, we will not be able to omit co-primality or squarefree-ness conditions in our computations; so we must choose terms a little carefully to obtain expressions that we can usefully work with.

When $k=o(\log\log x)$, we know from $\S 4$ that the $W=0$ term dominates the whole of the above sum. Roughly speaking, this is because $m$ cannot often have small prime factors, since numbers $a,b$ satisfying $\omega(a)=\omega(b)=k$ typically do not; and insisting that it should have certain large factors (when $W \geq 1$) greatly reduces the possibilities for $a,b$. For larger $k$ this reasoning fails, and e.g. the $W=1$ term becomes comparable with, and eventually much larger than, the $W=0$ term. We can write the $W=1$ term explicitly and simply, as
$$\sum_{p,q \leq x} \sum_{r < \min\{p,q\}, \atop r \textrm{ is prime}} \sum_{s < r, \atop s \textrm{ is prime}} 4 \#\{(pa'r,pa's) \in S_{p,k,x}^{2}\} \#\{(qc'r,qc's) \in S_{q,k,x}^{2}\}$$
$$= 4 \sum_{r \leq x, \atop r \textrm{ is prime}} \sum_{s < r, \atop s \textrm{ is prime}} \#\{t \in S_{k-1,x/r} : P(t)>r, \; r \nmid t, s \nmid t \}^{2},$$
recalling that $P(t)$ denotes the largest prime factor of $t$.

\vspace{12pt}
On a range $2 \leq k \leq R \log\log x$, for any fixed $R > 0$, the single pair $(r,s)=(3,2)$ gives a simple lower bound for the $W=1$ term. Recalling the notation $\overline{k}=k/L$ from the introduction, when $x$ is large enough (depending on $R$) we have $\overline{k} \leq 3R/2$. By Number Theory Result 5, then
\begin{eqnarray}
\#\{t \in S_{k-1,x/3} : P(t)>3, \; 2 \nmid t, 3 \nmid t \} & \gg_{R} & G(\overline{k-1}) \frac{(x/3) (\log\log(x/3))^{k-2}}{(k-2)! \log(x/3)} \nonumber \\
& \gg & \left(\frac{k-1}{\log\log x} \right) G(\overline{k-1}) \frac{x (\log\log x)^{k-1}}{(k-1)! \log x} \nonumber \\
& \gg & \overline{k} \#\{n \leq x : \omega(n)=\Omega(n)=k\}. \nonumber
\end{eqnarray}
In particular, if $\epsilon > 0$ and $\epsilon \log\log x \leq k \leq R \log\log x$, then for $x$ large enough the above is
$$ \gg_{\epsilon, R} \#\{n \leq x : \omega(n)=\Omega(n)=k\}. $$

When $k \leq (\log\log x)^{1.9}$ is larger, we must include more pairs of primes when seeking a lower bound. Qualitatively, numbers $t \leq x$ with many more than $\log\log x$ prime factors are almost always divisible by small primes like $2$ and $3$, so we need to look at a range of larger primes $r,s$. We will need a de Bruijn-type estimate, given in Chapter 7 of Montgomery and Vaughan~\cite{mv}: for any $\epsilon > 0$, and $y$ sufficiently large depending on $\epsilon$, one has
$$\#\{r \leq y : P(r) \leq \log^{2}y\} \leq y^{1/2 + \epsilon}.$$
Assume that $k > \log\log x$, so $\overline{k} > 1$. Combining the estimate with Number Theory Result 5, we find (with much to spare) that when $r \leq 4\overline{k}$, and $x$ is larger than an absolute constant,
\begin{eqnarray}
\#\{t \in S_{k-1,x/r} : P(t)>r, \; r \nmid t, s \nmid t \} & \geq & \#\{t \in S_{k-1,x/r} : r \nmid t, s \nmid t \} - \#\{t \leq x/r : P(t) \leq r\} \nonumber \\
& = & (1+o(1))\#\{t \in S_{k-1,x/r} : r \nmid t, s \nmid t \}. \nonumber
\end{eqnarray}
This implies that if a pair of primes $r,s$ satisfies $2\overline{k} \leq r \leq 4\overline{k}$, $\overline{k} < s < r$, then
\begin{eqnarray}
\#\{t \in S_{k-1,x/r} : P(t)>r, \; r \nmid t, s \nmid t \} & \gg & G\left(\overline{k} \left(1+O\left(\frac{1}{k} \right) \right) \right) \frac{x (\log\log x)^{k-2}}{(k-2)! \log x} \frac{1}{(r + \overline{k})} \nonumber \\
& \gg & G\left(\overline{k} \right) \frac{x (\log\log x)^{k-2}}{(k-2)! \log x} \frac{1}{(r + \overline{k})}. \nonumber
\end{eqnarray}
The second inequality uses the logarithmic derivative estimate at the end of $\S 3$.

We now sum over all such $(r,s)$, using the Chebychev-type estimate
$$0.9212y + O(\log y) \leq \sum_{p^{m} \leq y, \atop p \textrm{ prime}, m \in \N} \log p \leq 1.1056y + O(\log^{2}y) \; \; \; \textrm{if } y \geq 2,$$
and the fact that if $y > 1$ then there is a prime on the interval $(y, 2y)$. See Chapter 2 of Montgomery and Vaughan~\cite{mv}. It follows, as claimed, that the $W=1$ term is
$$ \gg G(\overline{k})^{2} \frac{x^{2} (\log\log x)^{2k-4}}{(k-2)!^{2} \log^{2}x} \sum_{2\overline{k} \leq r \leq 4\overline{k}, \atop r \textrm{ is prime}} \frac{r}{(r+\overline{k})^{2} \log r} \gg \left(\frac{\overline{k}}{\log(\overline{k}+2)} \#S_{k,x} \right)^{2}.$$

\vspace{12pt}
When $k \leq \delta \log x/(\log\log x)^{2}$ is even larger, we do not estimate $\#\{t \in S_{k-1,x/r} : P(t)>r, \; r \nmid t, s \nmid t \}$ so precisely. Instead, we rewrite the $W=1$ term as
$$4\sum_{t \in S_{k-1,x}} \sum_{u \in S_{k-1,x}} \sum_{r \leq x/t,x/u, \atop r \textrm{ is prime}} \textbf{1}_{r < P(t), P(u); \; r \nmid t,u} \sum_{s < r, \atop s \textrm{ is prime}} \textbf{1}_{s \nmid t,u}$$
$$ \geq 4 \sum_{t \in S_{k-1,x/3k\log k}} \sum_{u \in S_{k-1,x/3k\log k}} \textbf{1}_{P(t), P(u) > 3k\log k} \sum_{r \leq 3k\log k, \atop r \textrm{ is prime}} \textbf{1}_{r \nmid t,u} \sum_{s < r, \atop s \textrm{ is prime}} \textbf{1}_{s \nmid t,u}$$
$$ \geq k^{2} \left( \sum_{t \in S_{k-1,x/3k\log k}} \textbf{1}_{P(t) > 3k\log k} \right)^{2},$$
in view of the Chebychev-type lower bound for the sums over primes. (Here we noted that at most $2k-2$ of the primes $r$ less than $3k\log k$ are excluded by the presence of the indicator function.)

As before, but this time using Number Theory Result 3, we find
$$ \sum_{t \in S_{k-1,x/3k\log k}} \textbf{1}_{P(t) > 3k\log k} = (1+o(1))\#S_{k-1,x/3k\log k}. $$
Using both parts of Number Theory Result 4, and assuming as always that $x$ is large enough, this gives a lower bound for the $W=1$ term that is
$$\gg \frac{1}{\log^{2}k} (\#S_{k-1,x})^{2} \gg \left(\frac{k-1}{W \log k} \#S_{k,x} \right)^{2} \gg \left(\frac{\overline{k}}{\log k} \#S_{k,x} \right)^{2},$$
as asserted in Proposition 1.
\begin{flushright}
Q.E.D.
\end{flushright}

\section{Proof of Theorem 2}

\subsection{Strategy of the proof}
To establish Theorem 2, we will use a general fact about convergence in distribution of random variables: if $X_{n}, X$ are real valued random variables, and $X_{n} \stackrel{d}{\rightarrow} X$ as $n \rightarrow \infty$, then
$$\E g(X_{n}) \rightarrow \E g(X) \; \; \; \textrm{as } n \rightarrow \infty$$
whenever $g : \R \rightarrow \R$ is continuous and \emph{bounded}. This is sometimes given as the definition of convergence in distribution.

In particular, if $a \geq 0$ is fixed, and $X_{n} \stackrel{d}{\rightarrow} N(0,1)$ as $n \rightarrow \infty$, then
$$\E \min\{X_{n}^{2},a^{2}\} \rightarrow 1 - \frac{2a}{\sqrt{2\pi}} e^{-a^{2}/2} + 2(a^{2}-1)(1-\Phi(a)) \; \; \; \textrm{as } n \rightarrow \infty,$$
since $(2\pi)^{-1/2} \int_{a}^{\infty} y^{2} e^{-y^{2}/2} dy = (2\pi)^{-1/2}a e^{-a^{2}/2} + (1-\Phi(a))$.

With this in mind, it will suffice to find numbers $a \geq 0$ and $T_{a}$ such that $T_{a} > \frac{2a}{\sqrt{2\pi}} e^{-a^{2}/2} - 2(a^{2}-1)(1-\Phi(a))$, and
$$\E \min\{\widetilde{M}^{(k)}(x)^{2},a^{2}\} \leq 1-T_{a}$$
when $x$ is large (and $\epsilon \log\log x \leq k(x) \leq R \log\log x$). It may not be immediately clear that there is any sensible way to go about this, but general considerations at least suggest that the expression on the left should capture important information about $M^{(k)}(x)$. As remarked in $\S 2$, a normal approximation does hold for $\widetilde{M}^{(k)}(x)$ essentially when a sum of squared increments converges in probability to 1. This quantity is closely related to the square of $\widetilde{M}^{(k)}(x)$, e.g by Burkholder's inequality, as expounded in Hall and Heyde's book~\cite{hh}.

\vspace{12pt}
To proceed further, observe that if $q$ is any prime number then
\begin{eqnarray}
\E \min\{\widetilde{M}^{(k)}(x)^{2},a^{2}\} & = & \E \left( \E(\min\{\widetilde{M}^{(k)}(x)^{2},a^{2}\} \mid \epsilon_{2}, \epsilon_{3}, ..., \epsilon_{q}) \right) \nonumber \\
& \leq & \E \left( \min\{\E(\widetilde{M}^{(k)}(x)^{2} \mid \epsilon_{2}, \epsilon_{3}, ..., \epsilon_{q}), a^{2} \} \right) \nonumber \\
& = & 1 - \E \left( \max\{\E(\widetilde{M}^{(k)}(x)^{2} \mid \epsilon_{2}, \epsilon_{3}, ..., \epsilon_{q}) - a^{2}, 0 \} \right) \nonumber \\
& \leq & 1 - \p \left( \E(\widetilde{M}^{(k)}(x)^{2} \mid \epsilon_{2}, \epsilon_{3}, ..., \epsilon_{q}) \geq a^{2} + 1 \right). \nonumber
\end{eqnarray}
For given $x$, if we chose $q \geq x$ then the first inequality would be an equality. This suggests that computing $\E(\widetilde{M}^{(k)}(x)^{2} \mid \epsilon_{2}, \epsilon_{3}, ..., \epsilon_{q})$ as a function of $\epsilon_{2}, ..., \epsilon_{q}$ will be difficult when $q$ is large, (it amounts to explicitly determining the distribution of $\widetilde{M}^{(k)}(x)^{2}$), and we shall not attempt this.

However, for any given $q$, fixed before seeing what happens when $x \rightarrow \infty$, the computation becomes more feasible. In $\S 6.2$, we will carry this out to obtain an explicit, although not extremely enlightening, answer. Given this, there are two obvious ways to try to finish the proof of Theorem 2:
\begin{enumerate}
\item calculate the value (as a function of $x,k$) of $\E(\widetilde{M}^{(k)}(x)^{2} \mid \epsilon_{2}, \epsilon_{3}, ..., \epsilon_{q})$ for every possibility $\epsilon_{p}=\pm 1$, for some small values of $q$, and see if this leads to values $a, T_{a}$ with the desired properties;

\item choose values of $\epsilon_{p}$ that allow for good estimates of $\E(\widetilde{M}^{(k)}(x)^{2} \mid \epsilon_{2}, \epsilon_{3}, ..., \epsilon_{q})$ as $q$ becomes large, and then vary $q$ and $a$ until one can obtain suitable $T_{a}$.

\end{enumerate}
It is the second approach that leads to a proof of Theorem 2, as will be shown in $\S 6.2$ (just looking at $\E(\widetilde{M}^{(k)}(x)^{2} \mid \epsilon_{2}=1, \epsilon_{3}=1, ..., \epsilon_{q}=1)$). The first approach yields partial results, which are of some interest in that explicit numerical bounds for $\E \min\{\widetilde{M}^{(k)}(x)^{2},a^{2}\}$ are obtained, and we will discuss this briefly in $\S 6.3$.

\vspace{12pt}
The idea of using conditioning to explore the behaviour of $M^{(k)}(x)$ seems, to the author, rather natural, and it was already used in a heuristic way by Hough~\cite{ho}. He performs numerical simulations of the complete sum $M(x)$, both unconditionally and conditioning on a small number of $\epsilon_{p}$, and (looking at the empirical distribution functions) is led to conjecture that $M(x)$ ``...looks like a combination of conditional Gaussian distributions, whose variances depend on the value of $f$ on the first few primes.''

Hough~\cite{ho} also explains a heuristic of Chatterjee that $M(x)$ should not, in the limit, have a normal distribution: if it did, it seems likely that the distributions conditional on $\epsilon_{2}=1$ and $\epsilon_{2}=-1$ would also tend to normality, whilst these have distinct variances. The extent of the impact of $\epsilon_{p}$ on $M^{(k)}(x)$, when $p$ is `small' and $k$ is `not small', is what underlies both of these heuristics, and also the proof of Theorem 2. If the behaviour of $\widetilde{M}^{(k)}(x)^{2}$ were very close to, but not equal to, that of $N(0,1)^{2}$, our approach would probably not detect this.

\subsection{Completion of the proof}
We now implement the strategy proposed in $\S 6.1$, beginning by deriving a useable expression for $\E(\widetilde{M}^{(k)}(x)^{2} \mid \epsilon_{2}, \epsilon_{3}, ..., \epsilon_{q})$. The random multiplicative function $f$ will appear in this expression, but only applied to integers $N$ whose prime factors are at most $q$, as `short hand' for $\prod_{p|N} \epsilon_{p}$. We will also write
$$S_{\leq q, \leq k, x} := \{n \leq x: \omega(n)=\Omega(n) \leq k, P(n) \leq q\}.$$

Expanding the square, we see that $\E(M^{(k)}(x)^{2} | \epsilon_{2},...,\epsilon_{q})$ is
$$\sum_{n \in S_{k,x}} \sum_{m \in S_{k,x}} \E(f(n)f(m) | \epsilon_{2},...,\epsilon_{q})$$
$$ = \#S_{k,x} + \sum_{N \in S_{\leq q, \leq k,x}} f(N) \sum_{M \in S_{\leq q, \leq k,x}, \atop \omega(M)=\omega(N), M \neq N} f(M) \sum_{n \in S_{k-\omega(N), x/N} \cap S_{k-\omega(N), x/M}} \textbf{1}_{p \nmid n \forall p \leq q}. $$
Using Number Theory Result 5, we can replace the inner sum by
$$ (1+o(1)) G(\overline{k}(1+o(1))) \prod_{p \leq q} \left(1+ \frac{\overline{k}(1+o(1))}{p} \right)^{-1} \min\{x/N,x/M\} \frac{(\log\log x)^{k-1-\omega(N)}}{(k-1-\omega(N))! \log x}, $$
where the $o(1)$ terms are with respect to the limit process $x \rightarrow \infty$, for any fixed $q$. On the range of $k$ treated by Theorem 2, this is
$$ (1+o(1)) G(\overline{k}) \prod_{p \leq q} \left(1+ \frac{\overline{k}}{p} \right)^{-1} \min\{1/N,1/M\} \overline{k}^{\omega(N)} \frac{x(\log\log x)^{k-1}}{(k-1)! \log x};$$
so for any fixed $q$, we conclude that $\E(\widetilde{M}^{(k)}(x)^{2} \mid \epsilon_{2}, \epsilon_{3}, ..., \epsilon_{q})$ is
$$ 1 + 2 \prod_{p \leq q} \left(1+\frac{\overline{k}}{p} \right)^{-1} \sum_{N \in S_{\leq q, \leq k,x}} \frac{\overline{k}^{\omega(N)} f(N)}{N} \sum_{M \in S_{\leq q, \leq k,x}, \atop \omega(M)=\omega(N), M < N} f(M) + o(1).$$


\vspace{12pt}
Turning specifically to $\E(\widetilde{M}^{(k)}(x)^{2} \mid \epsilon_{2}=1, \epsilon_{3}=1, ..., \epsilon_{q}=1)$, it is clear that for any fixed $q$ this is at least
$$ 1 + 2 \prod_{p \leq q} \left(1+\frac{R}{p} \right)^{-1} \sum_{N \in S_{\leq q, \leq k,x}} \frac{\epsilon^{\omega(N)}}{N} \sum_{M \in S_{\leq q, \leq k,x}, \atop \omega(M)=\omega(N), M < N} 1 + o(1)$$
as $x \rightarrow \infty$. Moreover, provided $x \geq q^{2}-1$ is large enough (also depending on the value of $\epsilon$ in the statement of Theorem 2), this is at least
$$ 2 \prod_{p \leq q} \left(1+\frac{R}{p} \right)^{-1} \sum_{N \in S_{\leq q, 2,q^{2}-1}, N \geq 7} \frac{\epsilon^{2}}{N} \sum_{M \in S_{\leq q, 2, N-1}} 1, $$
which does not depend on $x$. We shall finish the proof of Theorem 2 by showing that, when $q$ is sufficiently large (depending on $\epsilon$ and $R$) and $x$ is sufficiently large, taking the preceding expression as $a^{2}+1$ yields
$$\frac{2a}{\sqrt{2\pi}} e^{-a^{2}/2} - 2(a^{2}-1)(1-\Phi(a)) < 2^{-q} < \p \left( \E(\widetilde{M}^{(k)}(x)^{2} \mid \epsilon_{2}, \epsilon_{3}, ..., \epsilon_{q}) \geq a^{2} + 1 \right).$$

Writing $p,q,r$ for prime numbers, observe that if $N \leq q^{2}$,
\begin{eqnarray}
\sum_{M \in S_{\leq q, 2, N-1}} 1 = \sum_{p \leq q} \sum_{p < r \leq \min\{q, (N-1)/p\}} 1 & = & \sum_{p < N/q} \sum_{p < r \leq q} 1 + \sum_{N/q \leq p \leq q} \sum_{p < r < N/p} 1 \nonumber \\
& = & \sum_{p < N/q} \sum_{p < r \leq q} 1 + \sum_{N/q \leq p \leq \sqrt{N}} \sum_{p < r < N/p} 1. \nonumber
\end{eqnarray}
One could estimate these sums very precisely using standard estimates for the distribution of prime numbers, but for us a fairly crude approach will suffice. If $2q < N \leq q^{2}/2$, then the first double sum is at least
$$ \sum_{p < N/q} \sum_{q/2 < r \leq q} 1 \gg \frac{q}{\log q} \sum_{p < N/q} 1 \gg \frac{N}{\log q \log(N/q)},$$
using the Chebychev-type lower bound quoted in $\S 5$. We then have
$$ \sum_{M \in S_{\leq q, 2, N-1}} 1 \gg \frac{N}{\log q \log(N/q + 2)} $$
whenever $7 \leq N \leq q^{2}$, since the left hand side is a non-decreasing function of $N$, and when $N \leq 2q$ the left hand side is simply $\#S_{2,N-1}$.

Using this bound, we find for $q$ sufficiently large that our choice of $a^{2}+1$ is
$$ \gg \prod_{p \leq q} \left(1+\frac{R}{p} \right)^{-1} \epsilon^{2} \frac{q^{2}}{\log^{4} q}. $$
Provided $q$ is larger than some absolute constant, a classical estimate of Mertens for $\sum_{p \leq y} 1/p$ reveals that the product is at least $(2 \log q)^{-R}$. See Chapter 2 of Montgomery and Vaughan's book~\cite{mv}, for example. This completes the proof of Theorem 2, with much to spare.
\begin{flushright}
Q.E.D.
\end{flushright}

\subsection{Computational aspects}
Recall that $\E(\widetilde{M}^{(k)}(x)^{2} \mid \epsilon_{2}, \epsilon_{3}, ..., \epsilon_{q})$ is
$$ 1 + 2 \prod_{p \leq q} \left(1+\frac{\overline{k}}{p} \right)^{-1} \sum_{N \in S_{\leq q, \leq k,x}} \frac{\overline{k}^{\omega(N)} f(N)}{N} \sum_{M \in S_{\leq q, \leq k,x}, \atop \omega(M)=\omega(N), M < N} f(M) + o(1),$$
for any fixed $q$. Because of the restriction that $\omega(M)=\omega(N)$ in the double sum, the number of terms that must be evaluated in this expression does not increase too quickly if $q$ is increased.

Assuming, for simplicity, that $\overline{k} \equiv 1$; that $x \geq \prod_{p \leq 29} p$ is large, so the sum over $N$ is over all squarefree numbers satisfying $P(N) \leq q$; and ignoring the $o(1)$ term (which can be made negligibly small, for computational purposes, by taking $x$ large enough); the author used Mathematica to evaluate the expression for values of $q$ up to $29$. This can be done at every point of the sample space, but in Table 1 we only present the maximal values attained.
\begin{table}[h]
\centering{
\begin{tabular}{l|c}
$q$ & $\E(\widetilde{M}^{(k)}(x)^{2} \mid \epsilon_{2}= ... = \epsilon_{q}=1)$ \\ \hline
2 & 1.000 \\
3 & 1.333 \\
5 & 1.806 \\
7 & 2.472 \\
11 & 3.249 \\
13 & 4.310 \\
17 & 5.603 \\
19 & 7.305 \\
23 & 9.378 \\
29 & 11.778 \\
\end{tabular}}
\caption{Conditional second moments, rounded to 3 decimal places.}
\end{table}

These calculations are enough to establish Theorem 2 when $\overline{k} \equiv 1$, and, since the expression computed is a continuous function of $\overline{k}$, also when $1-\delta \leq \overline{k} \leq 1 + \delta$ for some constant $\delta > 0$. We see
$$ \E \left( \max\{\E(\widetilde{M}^{(k)}(x)^{2} \mid \epsilon_{2}, \epsilon_{3}, ..., \epsilon_{29}) - 9, 0 \} \right) \geq \frac{2}{2^{10}} (11.777-9) \geq 0.0054;$$
and using tables of $\Phi(a)$, or using Mathematica,
$$\frac{6}{\sqrt{2\pi}} e^{-9/2} - 16(1-\Phi(3)) \leq 0.02660 - 0.02159 = 0.00501 < 0.0054.$$

\section{Sketch proof of Corollary 1}
First we show that Theorem 1 can be extended, by showing that $\widetilde{M}^{(k)}(x) - \widetilde{M}^{(\leq k)}(x)$ converges in probability to 0 as $x \rightarrow \infty$, when $k(x)=o(\log\log x)$. Introduce the temporary notation
$$\widehat{M}^{(\leq k)}(x) := M^{(\leq k)}(x)/\sqrt{\E M^{(k)}(x)^{2}}.$$
For any $\epsilon > 0$, we have
$$\p(|\widetilde{M}^{(k)}(x)-\widetilde{M}^{(\leq k)}(x)| \geq \epsilon) \leq \frac{2}{\epsilon^{2}} \left( \E(\widetilde{M}^{(k)}(x)-\widehat{M}^{(\leq k)}(x))^{2} + \E(\widehat{M}^{(\leq k)}(x)-\widetilde{M}^{(\leq k)}(x))^{2} \right),$$
by Chebychev's inequality. Using Number Theory Results 1 and 2, it is straightforward to show that the bracketed term is $o(1)$, as required.

The extension of Theorem 2 is based on the inequality
$$ \E(M^{(\leq k)}(x)^{2} | \epsilon_{2}=...=\epsilon_{q}=1) \geq \sum_{i=\min\{[k/2],[\log\log x/2]\}+1}^{\min\{k,[2\log\log x]\}} \E(M^{(i)}(x)^{2} | \epsilon_{2}=...=\epsilon_{q}=1).$$
If $k(x) \geq \epsilon \log\log x$, then each $i$ in the range of summation satisfies
$$ \min\{\epsilon/2,1/2\} \log\log x \leq i \leq 2\log\log x,$$
so as in $\S 6.2$ the summand is
$$\gg_{\epsilon} \left( \prod_{p \leq q} \left(1+\frac{2}{p} \right)^{-1} \sum_{N \in S_{\leq q, 2,q^{2}-1}, N \geq 7} \frac{1}{N} \sum_{M \in S_{\leq q, 2, N-1}} 1 \right) \#\{n \leq x : \omega(n)=\Omega(n)=i\},$$
provided $x$ is large enough. This suffices for the result.

\section{Sketch proof of Theorem 3}
We begin with the extension of Theorem 1 to the more general setting. We will write $g(n)$ for our generalised multiplicative function, where we allow non-Rademacher distributions for the underlying random variables $\epsilon_{p}$, whilst we continue to use $f(n)$ to denote a Rademacher random multiplicative function. Letting $p$ and $q$ be primes, we observe that
$$\E M_{p}^{(k)}(x)^{2} M_{q}^{(k)}(x)^{2} = \sum_{a \in S_{p,k,x}} \sum_{c \in S_{q,k,x}} \E g(a)^{2}g(c)^{2} + \sum_{a,b \in S_{p,k,x}, \atop a \neq b} \sum_{c,d \in S_{q,k,x}, \atop c \neq d} \E g(a)g(b)g(c)g(d). $$

Splitting the expectation in this way is just looking at $m=1$, and at other values of $m$, as in the proof of Theorem 1; we choose not to write it like this because some of the subsequent steps will be different. Summing the $m=1$ term over all pairs of primes $p,q$, using our assumptions about the $\epsilon_{p}$, we obtain
\begin{eqnarray}
\sum_{a \in S_{k,x}} \sum_{c \in S_{k,x}} \prod_{p | (a,c)} \E(\epsilon_{p}^{4}) & = & (\#S_{k,x})^{2} + O\left(\sum_{a \in S_{k,x}} \sum_{c \in S_{k,x}, (a,c) \neq 1} C^{\omega((a,c))} \right) \nonumber \\
& = & (\#S_{k,x})^{2} + O\left( \sum_{i=1}^{k} C^{i} \sum_{t \leq x, \omega(t)=i} (\#S_{k-i,x/t})^{2} \right). \nonumber
\end{eqnarray}
Similarly to Technical Lemma 1, one can show that (if $1 \leq i \leq k-1 \leq \log\log x$, and $x$ is large) the inner sum is
$$O\left(\frac{x^{2} (\log\log x + B)^{2k-2i-2}}{(k-i-1)!^{2}} \sum_{t \leq x, \omega(t)=i} \frac{1}{t^{2} \log^{2}(x/t+1)} \right) $$
$$ = O\left(\frac{x^{2} (\log\log x + B)^{2k-2i-2} (\log\log P_{i} + B)^{i-1}}{\log^{2}x (i-1)! (k-i-1)!^{2} P_{i} \log P_{i}} \right), $$
where $P_{i}$ is the least integer such that $\omega(P_{i})=i$, namely the product of the first $i$ primes. Thus the $m=1$ term is $(1+o(1))(\#S_{k,x})^{2}$, as $x \rightarrow \infty$ with $k(x)=o(\log\log x)$.

As in the proof of Theorem 1, it remains to show that the contribution from `other values of $m$', when summed over primes $p,q$, is $o((\#S_{k,x})^{2})$. Given $a,b,c,d \geq 2$, write $h(a,b,c,d)$ for the highest common factor of $a/P(a),b/P(b),c/P(c),d/P(d)$, and put $a'=a/h(a,b,c,d)$, and similarly $b',c',d'$. Then the contribution is at most
$$ \sum_{a,b \in S_{p,k,x}, \atop a \neq b} \sum_{c,d \in S_{q,k,x}, \atop c \neq d} C^{\omega(h(a,b,c,d))} \E g(a')g(b')g(c')g(d') $$
$$ \leq \sum_{i=0}^{k-1} C^{i} \sum_{t \leq x, \atop \omega(t)=i} \sum_{a,b \in S_{p,k-i,x/t}, \atop a \neq b} \sum_{c,d \in S_{q,k-i,x/t}, \atop c \neq d} \textbf{1}_{h(a,b,c,d)=1} \E g(a)g(b)g(c)g(d)$$
$$ \leq \sum_{i=0}^{k-1} C^{i+1} \sum_{t \leq x, \atop \omega(t)=i} \sum_{a,b \in S_{p,k-i,x/t}, \atop a \neq b} \sum_{c,d \in S_{q,k-i,x/t}, \atop c \neq d} \E f(a)f(b)f(c)f(d), $$
where the additional multiple of $C$ is needed on the third line because, sometimes, the smaller of $p$ and $q$ may divide all of $a,b,c,d$. Summing over $p,q$, and applying Lemma 1, one discovers that it would suffice if
$$ \sum_{i=1}^{k-W-1} \frac{C^{i} (\log\log x + B)^{2k-2i-2W-2} (\log\log P_{i} + B)^{i-1}}{(i-1)! (k-i-W-1)!^{2} P_{i} \log P_{i}} = O \left(\frac{(\log\log x + B)^{2k-2W-2}}{(k-W-1)!^{2}} \right) $$
and
$$ \sum_{i=1}^{k-W-1} \frac{C^{i} (2k-2i-2W-2)! (\log\log P_{i} + B)^{i-1}}{(i-1)! (k-i-W-1)!^{2} P_{i} \log P_{i}} = O \left(\frac{(2k-2W-2)!}{(k-W-1)!^{2}} \right), $$
uniformly for $1 \leq W \leq k-1$. Here the constants implicit in the ``big Oh'' notation will depend on the value of $C$. These results are straightforward to establish, in the manner of $\S 4.4$. 

\vspace{12pt}
The extension of Theorem 2 is less involved. For each prime $p$,
$$\p(\epsilon_{p} \geq 1/2) = \frac{1}{2} \E(\textbf{1}_{\epsilon_{p}^{2} \geq 1/4}) \geq \frac{(\E(\epsilon_{p}^{2} \textbf{1}_{\epsilon_{p}^{2} \geq 1/4}))^{2}}{2 \E \epsilon_{p}^{4}} \geq \frac{9}{32C}. $$
It follows as in $\S 6.2$ that, with probability at least $(9/32C)^{q}$, $\E(\widetilde{M}^{(k)}(x)^{2} \mid \epsilon_{2}, \epsilon_{3}, ..., \epsilon_{q})$ is greater than
$$ \prod_{p \leq q} \left(1+\frac{R}{p} \right)^{-1} \sum_{N \in S_{\leq q, \leq k,x}} \frac{\epsilon^{\omega(N)}}{N} \frac{1}{2^{\omega(N)}} \sum_{M \in S_{\leq q, \leq k,x}, \atop \omega(M)=\omega(N), M < N} \frac{1}{2^{\omega(N)}} \gg \frac{\epsilon^{2}}{2^{R}} \frac{q^{2}}{\log^{4+R}q} $$
for all large $x$ (depending on $q$, which itself must be larger than an absolute constant). The reader may check that this is still enough to establish the result, with much to spare.

\appendix
\section{Sketch proof of Technical Lemma 1}
We sketch the proofs of the four estimates making up Technical Lemma 1. To simplify the exposition, we will write
$$U:=[\log \sqrt{m}/\log 2], \;\; \textrm{ and } \;\; V:=[\log(M/m)/\log 2].$$

To deal with the double logarithms in a non-trivial way, we require the following bounds: if $n \in \N \cup \{0\}$, $C \geq 0$, and $x>1$, then
$$\int_{1}^{x} \frac{(\log t + C)^{n}}{t^{2}} dt \leq e^{C}n!;$$
and if $x > e$ then
$$\int_{e}^{x} \frac{(\log\log t + C)^{n}}{t \log^{2}t} dt \leq e^{C}n!.$$
These are easily established by making appropriate changes of variables, and induction on $n$.

For the first estimate, by Number Theory Result 1 the left hand side is
$$ O\left( \frac{(\log\log m + B)^{a-1}}{(a-1)! \log m} \sum_{i=0}^{U} \frac{(\log\log(2^{-i}\sqrt{m}+2) + C)^{n}}{\log^{2}(2^{-i-1}\sqrt{m}+1)} \right)$$
$$ = O\left( \frac{(\log\log m + B)^{a-1}}{(a-1)! \log m} \sum_{i=0}^{U} \frac{(\log\log(2^{U+1-i}+2) + C)^{n}}{\log^{2}(2^{U-i-1}+1)} \right).$$
Counting the terms of the sum backwards, we find it is
$$O\left( \int_{1}^{U+2} \frac{(\log(t+1) + C)^{n}}{t^{2}} dt \right),$$
whence the result.

By partial summation and Number Theory Result 1, the left hand side in the second estimate is
$$O\left(\frac{(\log\log(m+1)+C)^{n}(\log\log m + B)^{a-1}}{(a-1)! \log^{2}m} \right)$$
$$ + O\left(\int_{2}^{m} \frac{(\log\log(t+1)+C)^{n}}{t^{2} \log t} \cdot \frac{t(\log\log(t+1) + B)^{a-1}}{(a-1)! \log t} dt \right).$$  
Unless $m < 3$, when the estimate is trivial anyway, the first term may be omitted (at the cost of increasing the implicit constant, and values of $B$ and $C$ under the integral, by some fixed amounts). The result then follows as before.

The sum in the third estimate is at most
$$ O\left(\frac{M (\log\log M + B)^{a-1}}{(a-1)!} \sum_{i=0}^{V} \frac{2^{i+1}m}{M \log(2^{i+1}m) \log^{2}(1+2^{i})} \right).$$
Removing a factor $1/\log M \log^{2}(1+M/m)$ from the sum, we are left to bound
$$ \sum_{i=0}^{V} 2^{i-V} \frac{\log(2^{V+1}m)}{\log(2^{i+1}m)} \frac{\log^{2}(1+2^{V})}{\log^{2}(1+2^{i})} = O\left(\sum_{i=0}^{V} 2^{i-V} \left( \frac{V+1}{i+1} \right)^{3} \right).$$
It is easy to see that this quantity can be bounded independently of $m$ and $M$, e.g. by considering the ratio of consecutive summands.

For the fourth estimate, it will suffice to show that
$$\sum_{i=0}^{[\log m/\log 2]} \frac{\log m}{2^{[\log m/\log 2]-i} (i+1)} \left(\frac{\log\log(N/2^{i})+C}{\log\log(N/m)+C}\right)^{t}$$
has a bound depending only on $C$ and $D$. However, the ratio of the $i$ and $i-1$ summands is at least
\begin{eqnarray}
\frac{2i}{i+1} \left(1+\frac{\log(1-\log 2/\log(N/2^{i-1}))}{\log\log(N/2^{i-1})+C}\right)^{t} & \geq & \frac{2i}{i+1} \left(1+O\left(\frac{1}{\log(N/2^{i-1}) \log\log(N/m)} \right) \right)^{t} \nonumber \\
& \geq & \frac{2i}{i+1} \left(1+O\left(\frac{1}{\log(N/m) \log\log(N/m)} \right) \right)^{t}. \nonumber
\end{eqnarray}
This is more than $5/4$, say, provided that $i \geq 2$ and $N/m$ is larger than a constant depending on $D$ only. If $N/m$ is not so large, we have
$$ \log\log(N/2^{i}) \ll_{D} \log\log(m/2^{i} + 3) \;\;\;\; \textrm{and} \;\;\;\; t \ll_{D} 1.$$
The reader can check that the sum over $i$ is dominated by $O_{D}(1)$ terms with $i$ largest, and the result still follows.

\vspace{12pt}
\noindent {\em Acknowledgements.} The author would like to thank his PhD supervisor, Ben Green, for suggesting that he work on this problem and for many discussions about it.

\end{document}